\documentclass[10pt]{article}
\usepackage{amsmath}
\usepackage{amsthm}

\newtheorem{theorem}{Theorem}[section]
\newtheorem{lemma}[theorem]{Lemma}
\newtheorem{proposition}[theorem]{Proposition}

\newcommand{\Endo}{\mathop{\mathrm {End}}\nolimits}
\newcommand{\LNalt}{\mathop{\mathrm {LN_{alt}}}\nolimits}
\newcommand{\InnDer}{\mathop{\mathrm{InnDer}}\nolimits}
\newcommand{\spann}{\mathop{\mathrm {span}}}

\newcommand{\Lie}{\mathop{\mathcal{L}}\nolimits}
\newcommand{\modulus}[1]{\enskip (\mathop{\mathrm {mod}} #1)}
\newcommand{\id}{\mathop{\mathrm{Id}}\nolimits}
\newcommand{\ad}{\mathop{\mathrm{ad}}\nolimits}
\newcommand{\trace}{\mathop{\mathrm{tr}}\nolimits}

\begin{document}

\title{Right ideals in non--associative universal enveloping algebras of Lie triple systems}
\author{Jos\'e M. P\'erez-Izquierdo\footnote{Supported by the CONACyT (U44100-F)}}
\date{}
\maketitle

\begin{abstract}
We prove that the only proper right ideal of the universal enveloping algebra of a
finite--dimensional central simple Lie triple system over a field of characteristic zero is its
augmentation ideal.
\end{abstract}

\section{Introduction}

\bigskip
\emph{In this paper we will assume that the base field is algebraically closed of characteristic zero.}
\bigskip

Any Lie algebra with the trilinear product $[x,y,z] = [[x,y],z]$ can
be considered as a Lie triple system. A Lie triple system
(abbreviated as L.t.s.) is a vector space $T$ equipped with a
trilinear product $[x,y,z]$ satisfying the followings identities:
\begin{eqnarray}
&[x,x,y] = 0,&\\ &[x,y,z] + [y,z,x] + [z,x,y] = 0,\label{Jacobi.eq}&\\ &[a,b[x,y,z]] =
[[a,b,x],y,z] + [x,[a,b,y],z] + [x,y,[a,b,z]].\label{derivation.eq}&
\end{eqnarray}
The first and the second identities are reminiscent of Lie algebras.
The third identity says that the maps $D_{a,b}\colon x\mapsto
[a,b,x]$ are derivations of $T$. These maps are called inner
derivations and they form the Lie algebra of inner derivations
$\mathrm{InnDer}(T)$. Roughly, Lie triple systems may be thought of
as subspaces of Lie algebras stable under the product $[[x,y],z]$.
This is due to the fact that for any Lie triple system $T$ the
vector space $L(T) = \InnDer(T)\oplus T$ is a Lie algebra with the
product defined by $[a,b] = D_{a,b}$ and $[D_{a,b},c] =
[[a,b],c]=[a,b,c]$ for any $a,b,c\in T$, and by the condition that
$\InnDer(T)$ a subalgebra. The map $\sigma \colon D_{a,b} + c\mapsto
D_{a,b} -c$ is an automorphism of $L(T)$ with $\sigma^2 = \id$.
Thus, Lie triple systems are often considered as the negative
eigenspaces of involutions of Lie algebras. The classification of
these involutions has led to the classification of semisimple
L.t.s.\ by Lister \cite{L}. Later, Faulkner obtained Lister's
classification of  simple L.t.s.\ by means of Dynkin diagrams
\cite{F}.

Given any (unital) nonassociative algebra $A$ the generalized left
alternative nucleus
\begin{displaymath}
\LNalt(A) = \{ a\in A\,\vert\, (a,x,y) = -(x,a,y)\,\forall x,y \in A\},
\end{displaymath}
where $(x,y,z) = (xy)z-x(yz)$, is a L.t.s.\  with the product
$[a,b,c] = a(bc) - b(ac) -c(ab) + c(ba)$. In analogy with the usual
universal enveloping algebras of Lie algebras, for any L.t.s.\  $T$
there exists a unital algebra $U(T)$ and a monomorphism of L.t.s.\
$\iota \colon T\rightarrow \LNalt(U(T))$ with the additional
property $\iota(a)\iota(b) -\iota(b)\iota(a) = 0$ for any $a,b\in
T$. The pair $(U(T),\iota)$ is universal with respect to
homomorphisms of L.t.s.\ $T\rightarrow \LNalt(A)$ with this
property. To simplify the notation we will write $a$ instead of
$\iota(a)$. Thus,
\begin{displaymath}
[a,b,c] = a(bc) - b(ac) = a(bc) -(ab)c + (ba)c-b(ac) =
-2(a,b,c)
\end{displaymath}
in $U(T)$.

The universal enveloping algebra $U(T)$ is a nonassociative Hopf
algebra or H--bialgebra (see \cite{P,SU}) in the sense that it is a
nonassociative bialgebra with a left and a right division
$x\backslash y$ and $x/y$ defined by
\begin{eqnarray*}
&\sum x_{(1)}\backslash x_{(2)}y = \epsilon (x) y = \sum x_{(1)}\cdot x_{(2)}\backslash y, &\\
&\sum yx_{(1)}/ x_{(2)} = \epsilon (x) y = \sum y/x_{(1)}\cdot x_{(2)},&
\end{eqnarray*}
where $\epsilon$ denotes the counit and $\sum x_{(1)}\otimes
x_{(2)}$ stands for the image of $x$ under the comultiplication
\cite{Sw}. $U(T)$ is a coassociative, cocommutative bialgebra, the
subspace of primitive elements being $T$, and it admits a
Poincar\'e--Birkhoff--Witt type basis. Associativity is, however,
replaced by the weaker identity
\begin{displaymath}
\sum x_{(1)}(y\cdot x_{(2)}z) = \sum x_{(1)}(y x_{(2)})\cdot z.
\end{displaymath}
The graded algebra $\mathrm{Gr}(T)$ associated to the coradical
filtration
\begin{displaymath} U(T) = \bigcup_{n=0}^\infty U(T)_n,
\end{displaymath} where $U(T)_n$ is the linear span of all possible
products of at most $n$ elements in $T$, is isomorphic to the
symmetric algebra $S(T)$ on $T$. A map similar to the antipode of
universal enveloping algebras of Lie algebras is also available since the
automorphism $a\mapsto -a$ of any L.t.s.\  $T$ induces an
automorphism $x\mapsto S(x)$ of $U(T)$ of order $2$. The left and
right division are written as
\begin{displaymath}
x\backslash y = S(x) y\quad \mathrm{ and }\quad y/x = \sum S(x_{(3)})\cdot (x_{(1)}y) S(x_{(2)}).
\end{displaymath}

For any non-trivial L.t.s.  $T$  its universal enveloping algebra $U(T)$ is
infinite--dimensional. For finite--dimensional unital algebras $A$
the existence of embeddings $\iota\colon T\rightarrow \LNalt(A)$
with the property $\iota(a)\iota(b)-\iota(b)\iota(a) = 0 \, \forall
a,b\in T$, (Ado's Theorem) was studied in \cite{MP}: it turns out
that only finite--dimensional nilpotent L.t.s.\ may admit such
embeddings. The lack of such embeddings is equivalent to the
non-existence of ideals of finite codimension in $U(T)$ with trivial
intersection with $T$. This result motivated in \cite{MP} the
following conjecture, verified in several cases by direct
calculations in Poincar\'e--Birkhoff--Witt bases:

\bigskip

\noindent\textbf{Conjecture.}\enskip\emph{The only proper ideal of the universal enveloping
algebra of a  simple Lie triple system is its augmentation ideal.}

\bigskip

The goal of this paper is to prove this conjecture by establishing
the following stronger result:

\begin{theorem}\label{main.th}
The only proper right ideal of the universal enveloping algebra of a
finite--dimensional  simple Lie triple system over a field of
characteristic zero is its augmentation ideal.
\end{theorem}

The behavior of the left ideals is, however, rather different. The
Harish-Chandra isomorphism for symmetric spaces, recast in our
setting, implies that the right associative nucleus of $U(T)$, each
of whose elements determines a left ideal, is isomorphic as an
algebra to the algebra of invariants of the restricted Weyl group
of the symmetric Lie algebra $(L(T),\sigma)$. 
By the Theorem of Kostant and Rallis on separation of variables for
isotropy representations, $U(T)$ is a free right module for the
right associative nucleus. This subject will be discussed in detail
elsewhere.

\bigskip

\noindent{\textbf{Acknowledgment.} I would like to express my gratitude to Jacob Mostovoy for careful reading of the manuscript.}

\section{The proof}

 In the course of the proof we use some results whose proofs are postponed to later sections for
clarity.

Recall that for any $c\in T$ the subalgebra generated by $c$ is
associative, so $c^n$ is a well defined element of $U(T)$. In fact,
$(c^i, c^j,x) = 0$ for any $x \in U(T)$ \cite{MP}.

\begin{lemma}
\label{derivation.lem}
For any $n\geq 0$ and $a,b,c\in T$, we have
\begin{displaymath}
(c^n,a,b) \equiv n c^{n-1}(c,a,b) \modulus{U(T)_{n-2}}.
\end{displaymath}
\end{lemma}

The meaning of this lemma is better understood by defining $
R_{a,b}\colon T \rightarrow T$  $c \mapsto [c,a,b]$
or its extension to the whole $U(T)$
\begin{displaymath}
R_{a,b}(x) = - 2(x,a,b).
\end{displaymath}
These maps preserve the filtration of $U(T)$, so they induce
corresponding maps $R_{a,b}^{\mathrm{gr}}$ of zero degree on the
graded algebra ${\mathrm{Gr}}(U(T)) = \bigoplus_{n=0}^{\infty}
U(T)_n/U(T)_{n-1}$, which is isomorphic to $S(T)$. Since the powers
$\{c^n\,\vert\, c\in T, n\geq 0\}$ span $S(T)$, Lemma
\ref{derivation.lem} says that $R_{a,b}^{\mathrm{gr}}$ is the
derivation on $S(T)$ induced by the map $c\mapsto [c,a,b]$ on $T$.

The unifying feature of  simple L.t.s.\  that we use in our proof is
recorded in the following theorem: 

\begin{theorem}
\label{Endo.th} Let $T$ be a finite--dimensional  simple L.t.s.\ over
a field of characteristic zero. Then $\Endo(T)$ is the Lie algebra
generated by the maps $R_{a,b}\colon c\mapsto [c,a,b]$ $(a,b\in T)$.
\end{theorem}

Any finite--dimensional simple L.t.s.\  contains a two--dimensional
subsystem $S_2 =\spann\langle e,f\rangle$ with the product given by
\begin{displaymath}
[e,f,e] = 2e \quad \mathrm{ and } \quad [e,f,f] = -2f.
\end{displaymath}
We can reduce the proof of Theorem~\ref{main.th} to the case $T=
S_2$ as follows. Any nonzero right ideal $I$ of $U(T)$ is filtered
by $I_n = I\cap U(T)_n$ so it gives rise to a (two--sided) graded
ideal
\begin{displaymath} I^{\mathrm{gr}}=\bigoplus_{n=0}^\infty (I_n + U(T)_{n-1})/U(T)_{n-1}\cong
\bigoplus_{n=0}^\infty I_n/I_{n-1}
\end{displaymath}
that is stable under the derivations $R^{\mathrm gr}_{a,b}$ $\forall
a,b\in T$. By Theorem \ref{Endo.th}, the action of these operators
on $\mathrm{Gr}(T)$ can be identified with the natural action of
$\Endo(T)$ on $S(T)$. This module decomposes as direct sum of
irreducible modules $S(T)\cong \oplus_{n=0}^\infty V(n \lambda_1)$
where $V(n\lambda_1)$ corresponds to the homogeneous polynomials of
degree $n$. Since $I^{\mathrm gr}$ is a submodule as well as a
graded ideal then there must exist $N$ such that
\begin{displaymath}
I^{\mathrm gr} = \bigoplus_{n=N}^{\infty} U(T)_n/U(T)_{n-1}.
\end{displaymath}
In terms of $I$ we get the decomposition
\begin{equation}
\label{ideal.eq}
U(T) = U(T)_{N-1}\oplus I
\end{equation}
which shows that the codimension of $I$ is finite. Therefore, given
a subsystem $S_2=\spann\langle e,f\rangle \subseteq T$, the right
ideal $I\cap U(S_2)\subseteq U(S_2)\subseteq U(T)$ has finite
codimension in $U(S_2)$. Assuming the truth of Theorem~\ref{main.th}
in the case of $S_2$, we get that $I\cap T \neq 0$. Since $T$ is
simple this intersection generates all $T$ under the action of the
operators $R_{a,b}$. Thus, $I $ is either the augmentation ideal or
the whole $U(T)$ as desired.

\section{Proof of Lemma \ref{derivation.lem}}

The main tool in computing products in $U(T)$ is the operator identity
\begin{equation}
\label{Jordan.eq}
 L_{ax + xa} = L_aL_x + L_x L_a
\end{equation}
for all $a\in T$ and $x\in U(T)$ \cite{MP}. To use it we first observe that
\begin{eqnarray*}
ac^n &=& ac\cdot c^{n-1} - (a,c,c^{n-1}) = ca\cdot c^{n-1} -(a,c,c^{n-1})\\ &=& c\cdot ac^{n-1} +
(c,a,c^{n-1}) -(a,c,c^{n-1}) = c\cdot ac^{n-1} - 2(a,c,c^{n-1})\\ &=& -2(a,c,c^{n-1}) - 2c\cdot
(a,c,c^{n-2}) -\cdots - 2 c^{n-2}(a,c,c) + c^n a,
\end{eqnarray*}
so
\begin{displaymath}
c^n a =\frac{1}{2}(ac^n + c^n a) + \sum_{i=0}^{n-2} c^i(a,c,c^{n-1-i}).
\end{displaymath}
The associator $(c^n,a,b)$ is then written as
\begin{eqnarray*}
(c^n,a,b) &=& \frac{1}{2} L_{ac^n + c^n a}(b)+ \sum_{i=0}^{n-2} L_{c^i (a,c,c^{n-1-i})}(b) -
L_{c^n}L_a(b)
\end{eqnarray*}
that by (\ref{Jordan.eq}) gives
\begin{eqnarray*}
(c^n,a,b)&=& \frac{1}{2}[L_a,L_{c^n}](b) + \sum_{i=0}^{n-2}L_{c^i
(a,c,c^{n-1-i})}(b)\\ &=& \frac{1}{2} \sum_{i=0}^{n-1}L^i_c [L_a,L_c]L^{n-1-i}_c(b) +
\sum_{i=0}^{n-2}L_{c^{i}(a,c,c^{n-1-i})}(b)
\end{eqnarray*}
where we have used that $L_{c^n} = L^n_c$. Since $D_{a,b}=[L_a,L_b]$ is a derivation of
$U(T)$ \cite{MP} we get
\begin{eqnarray}
 \label{associator.eq}
 (c^n,a,b)
&=& \frac{1}{2} \sum_{i=0}^{n-1}c^i D_{a,c}(c^{n-1-i}b) + \sum_{i=0}^{n-2}c^i(a,c,c^{n-1-i})\cdot
b\nonumber\\ &=& \frac{1}{2} \left( \sum_{i=0}^{n-1} c^i \cdot D_{a,c}(c^{n-1-i})b + c^{i}\cdot
c^{n-1-i}[a,c,b]\right) \nonumber\\ && + \sum_{i=0}^{n-2}c^i(a,c,c^{n-1-i})\cdot b\nonumber\\ &=&
\frac{n}{2}c^{n-1}[a,c,b] + \frac{1}{2}\sum_{i=0}^{n-2} c^i\cdot D_{a,c}(d^{n-1-i})b \nonumber\\
&& - \frac{1}{2}\sum_{i=0}^{n-2} c^i D_{a,c}(c^{n-1-i}) \cdot b\nonumber \\ &=&
\frac{n}{2}c^{n-1}[a,c,b] - \frac{1}{2}\sum_{i=0}^{n-2}(c^i,D_{a,c}(c^{n-1-i}),b).
\end{eqnarray}
The derivation $D_{a,c}$ preserves the filtration of $U(T)$ so
$(c^i,D_{a,c}(c^{n-1-i}),b) \in U(T)_{n-2}$ and $(c^n,a,b) \equiv
\frac{n}{2} c^{n-1}[a,c,b] = -n c^{n-1}(a,c,b) = nc^{n-1}(c,a,b)$
modulo $U(T)_{n-2}$.

\section{The two--dimensional case}

The case $T= S_2$ is simple and illustrative because we can perform
computations in the Poincar\'e--Birkhoff--Witt basis
$\{e^if^j\,\vert\, i,j\geq 0\}$. In this case Theorem~\ref{Endo.th}
encodes the following correspondence between maps and coordinate
matrices in the basis $\{ e, f\}$:
  \begin{displaymath}
 R_{e,e}\equiv \begin{pmatrix}0 & -2 \\0 & 0\end{pmatrix},  R_{e,f} \equiv \begin{pmatrix}0 & 0
\\ 0& 2\end{pmatrix},R_{f,e}\equiv \begin{pmatrix}2 &0 \\ 0 & 0\end{pmatrix} \mathrm{\  and\  } R_{f,f}
\equiv
\begin{pmatrix}0 & 0 \\ -2 & 0\end{pmatrix}.
\end{displaymath}
Fix $I$ a nonzero right ideal of $U(S_2)$. By (\ref{ideal.eq}), $I$
contains an element of the form $e^N+ \alpha_1 e^{N-1}+ \cdots +
\alpha_{N-1} e + \alpha_N$. Moreover, since (\ref{associator.eq})
implies that $R_{f,e}(e^n)=-2(e^n,f,e) = 2 n e^n$ then $I$ must
contain the power $e^N$, although $e^{N-1}\not\in I$. The result is
a consequence of the following proposition which shows that in this
situation $N = 1$ or $N = 0$, so $I$ is either the augmentation
ideal or the whole $U(S_2)$.

\begin{proposition}
In $U(S_2)$ we have that $(e^n,f,f)e = n e^n f - n(n-1)e^{n-1}$ $\forall n\geq 0$.
\end{proposition}
\begin{proof}
The cases $n = 0$ and $n=1$ are obvious, so we assume that $n\geq
2$. On one hand, $(e^n,f,f)e$ is an eigenvector  of $D_{e,f}$ with
eigenvalue $2n -2$, so it is a linear combination of $\{ e^{n-1 +
i}f^i\,\vert\, i\geq 0\}$. Since it also belongs to $U(S_2)_{n+1}$,
we deduce that $(e^n,f,f)e = \alpha e^n f + \beta e^{n-1}$ for some
$\alpha,\beta\in F$. On the other hand, by (\ref{associator.eq})
$(e^n,f,f) e = n e^{n-1}f\cdot e = n e^n f - n [e,e^{n-1}f]$ and, by
Lemma 24 in \cite{MP}, this is congruent to $n e^n f - n(n-1)
e^{n-1}$ modulo $U(S_2)_{n-3}$. Therefore, $\alpha = n$ and $\beta =
-n(n-1)$.
\end{proof}

\section{Proof of Theorem \ref{Endo.th}}

Let $T$ be a finite--dimensional simple L.t.s.\ and $L=L(T) =
\InnDer(T)\oplus T$ --- its standard embedding Lie algebra. The
results of the previous section allow us to assume that $\dim{T}\geq
3$. By \cite{L} $L(T)$ is a semisimple Lie algebra so its Killing
form $K(\,,\,)$ is nondegenerate. The corresponding automorphism
$\sigma$ of order $2$ preserves the Killing form so $\InnDer(T)$ and
$T$ are orthogonal and the restriction of $K(\,,\,)$ to these
subspaces is nondegenerate. The adjoint of $R_{a,b} =
\ad_b\ad_a\vert_T$ relative to $K(\,,\,)$ is $R_{b,a}$ and
$$2\trace(R_{a,b}) = 2\trace(\ad_b\ad_a\vert_T) = K(a,b).$$

As an $\InnDer(T)$--module, $T$ is either irreducible or it
decomposes into a direct sum of two irreducible submodules $T = T_1
\oplus T_2$ with the further restriction that $[T_1,T_1] = 0 =
[T_2,T_2]$. In the latter case the center $Z(\InnDer(T))$ of
$\InnDer(T)$ is one--dimensional and it is spanned by the map acting
as the identity $\id$ on $T_1$ and as $-\id$ on $T_2$. The
skew--symmetry of this map with respect to the Killing form implies
that $T_1$ and $T_2$ are isotropic relative to the Killing form, so
$T_2\cong T^*_1$, the dual module of $T_1$.

One important property of simple L.t.s.\  that we will use is that
for any $0\neq x_{\eta} \in T_\eta$ with $\eta$ a nonzero weight of
$T$ relative to some Cartan subalgebra of $\InnDer(T)$ there exists
$x_{-\eta} \in T_{-\eta}$ such that $\spann\langle
x_{\eta},x_{-\eta}, [x_{\eta},x_{-\eta}]\rangle$ is a
three--dimensional Lie algebra isomorphic to $sl(2)$ \cite{F}. In
particular, $[x_{\eta},x_{-\eta},x_{\eta}]\neq 0$.

Let $\Lie$ be the Lie algebra generated by $\{ R_{a,b}\,\vert\,
a,b\in T\}$.  We will use the notation
\begin{displaymath}
\tau_{x,y}\colon z \mapsto K(y,z) x, \quad
\sigma_{x,y} = \tau_{x,y} + \tau_{y,x}, \textrm{ and }
\lambda_{x,y} =  \tau_{x,y} - \tau_{y,x}.
\end{displaymath}

Assume that $T$ is a $\InnDer(T)$--irreducible module. Fix $\mu$ to
be the highest weight of $T$ relative to a basis of a root system of
$\InnDer(T)$ and $0\neq x_\mu \in T_\mu$. For any weight $\eta$ of
$T$ and $x_\eta \in T_\eta$,  $[x_\eta,x_\mu,x_\mu] \in T_{2\mu +
\eta}$ so  it vanishes if $\eta \neq -\mu$. Since in case that $\eta
= -\mu$ we have already noticed that $[x_{-\mu},x_\mu,x_\mu] \neq
0$, then $R_{x_\mu,x_\mu}$ is a nonzero multiple of
$\tau_{x_{\mu},x_{\mu}}$. Therefore, $\tau_{x_{\mu},x_{\mu}}$ as
well as  $\tau_{x_{-\mu},x_{-\mu}}$ and  $
\lambda_{x_{\mu},x_{-\mu}}
=[\tau_{x_{\mu},x_{\mu}},\tau_{x_{-\mu},x_{-\mu}}] $ belong to
$\Lie$.

For any map $d \in so(T)$ skew--symmetric with respect to the
Killing form we have that $[d,\tau_{x,y}] = \tau_{d(x),y} +
\tau_{x,d(y)}$ so by using root vectors of positive roots, and due
to the irreducibility of $T$, we easily get that $\lambda_{x_\mu,y}
\in \Lie$ for any $y \in T$. From this we can obtain that $so(T)
=\spann\langle \lambda_{x,y}\,\vert\, x,y\in T\rangle \subseteq
\Lie$.

The traceless map $0\neq R_{x_\mu,x_\mu}\in \Lie$ is symmetric with
respect to the Killing form. By the usual decomposition of $sl(T)$
into the direct sum of two irreducible submodules (symmetric and
skew--symmetric maps) with respect to the  action of $so(T)$ by
commutation it follows that $sl(T) \subseteq \Lie$. Since $2
\trace(R_{x,x}) = K(x,x)$ then $\Lie$ also contains the identity,
so $\Lie = \Endo(T)$.

Let us assume now that $T$ is not irreducible. We have $$\Endo(T) =
\Endo(T)_0\oplus \Endo(T)_1$$ with $$\Endo(R)_0 = \{ f\in
\Endo(T)\,\vert\, f(T_j) \subseteq T_{j}, j =1,2\}$$ (even maps) and
$$\Endo(R)_1 = \{ f\in \Endo(T)\,\vert\, f(T_j) \subseteq T_{3-j}, j
=1,2\}$$ (odd maps).

Fix $\mu$ to be the highest weight of $T_1$ with respect to a basis
of a root system of $\InnDer(T)$, $0\neq x_\mu \in (T_1)_\mu$ and
$0\neq x_{-\mu}\in (T_2)_{-\mu}$ (recall that $T_2 \cong T^*_1$).
Since $[T_1,T_1] =0$ and $T_2 \cong (T_1)^*$ then
$[x_\eta,x_\mu,x_\mu] = 0$ for any weight $\eta \neq -\mu$. As
before, this implies that $\lambda_{x_{\mu},x_{-\mu}} \in \Lie$
which ultimately gives that $$so(T)_0 := so(T)\cap \Endo(T)_0 =
\lambda_{T_1,T_2} \subseteq \Lie.$$   The dimension of $T_1$ is
$\geq 2$, so  under commutation $\Endo(T)_0$ is the sum of three
irreducible $so(T)_0$--modules: $so(T)_0$, the traceless even
symmetric maps and the scalar maps. As above $\{ R_{x,y} +
R_{y,x}\,\vert\, x \in T_1, y\in T_2\}$ is not contained inside
$sl(T)$, so $\id \in \Lie$. In order to prove that $\Endo_0(T)
\subseteq \Lie$ we only have to check that $\Lie_0 \neq so(T)_0
\oplus F \id$. By the contrary, if $R_{x,y}+ R_{y,x}\subseteq F \id$
for any $x\in T_1, y \in T_2$ then taking traces we obtain that
$R_{x,y} + R_{y,x} = \frac{K(x,y)}{\dim T_1} \id$.  Given $x' \in
T_1$, by (\ref{Jacobi.eq})
\begin{displaymath}
\frac{K(x,y)}{\dim T_1} x' = [x',x,y] + [x',y,x] = [x',y,x]  = [x,y,x'] = \frac{K(x',y)}{\dim T_1} x
\end{displaymath}
which contradicts that $\dim T_1 \geq 2$. Therefore,
$\Endo_0(T)\subseteq \Lie$.

$\Endo_1(T)$ decomposes as the direct sum of two irreducible
submodules under the action of $\Endo_0(T)$ by commutation, namely
those maps that kill $T_1$ and those that kill $T_2$. Since
$R_{x_{\mu},x_{\mu}}$ and $R_{x_{-\mu},x_{-\mu}}$ are nonzero
elements of these types then $\Endo_1(T) \subseteq \Lie$, so $\Lie =
\Endo(T)$ as desired.

\end{document}